
\input amstex
\documentstyle{amsppt}
\magnification = \magstep 1
\NoBlackBoxes
\NoRunningHeads

\vsize=8.9truein
\hsize=6.4 true in
\hcorrection{ -.1875 true in}
\vcorrection{  -.25 true in}

\define\({\left (}
\define\){\right )}

\define\brq{^{[q]}}
\define\brt{^{[t]}}
\define\q0{{q_0}}
\redefine\vec#1 #2{#1_1,\ldots, #1_{#2}}

\define\bdf{\bold f}

\define\bz{\bold z}

\define\m{\bold m}
\define\n{\bold n}

\define\Q{\Bbb Q}
\define\Z{\Bbb Z}

\define\Rq{R^{1/q}}
\define\Sq{S^{1/q}}

\define\inc{\subseteq}
\define\1{^{-1}}
\define\8{{\infty}}
\define\Ro{R^{\circ}}
\define\defn#1{{\sl #1}}

\define\hatR{{\widehat R}}
\define\Rhat{{\widehat R}}

\define\al{\alpha}
\define\eps{\epsilon}

\def\ord {\operatorname{ord}\nolimits}

\topmatter
\title Extension of weakly and strongly $F$-regular rings
by flat maps 
 \endtitle
\author  Ian M.~Aberbach \endauthor
\address Department of Mathematics, University of Missouri,
Columbia, MO  65211 \endaddress
\email aberbach\@math.missouri.edu \endemail
\thanks  The  author was partially supported by the NSF.\endthanks
\endtopmatter
\document

\head \S 1.  Introduction \endhead

Throughout this paper all rings will be Noetherian of positive
characteristic $p$.  Hence tight closure theory \cite{HH1--4} takes
a prominent place (see \S 2 for tight closure definitions and terminology).
  The purpose of this note is to help answer the following
question:  if $R$ is weakly (resp.~strongly) $F$-regular  and
$\phi:R \to S$ is a flat map then under
what conditions on the fibers is $S$ weakly (resp.~strongly) $F$-regular.
This question (among many others) is raised in \cite{HH4} in section 7.
It is shown there that if $\phi$ is a flat map of local rings, $S$
is excellent and the generic and closed fibers are regular then weak
$F$-regularity of $R$ implies that of $S$ (Theorem 7.24).  One of our
main results weakens the hypotheses considerably.

\proclaim{Theorem 3.4}  Let $\phi:(R,\m) \to (S,\n)$ be a flat map.  Assume that
$S/\m S$ is Gorenstein and $R$ is weakly $F$-regular and 
Cohen-Macaulay.  Suppose that either
\roster
\item  $c \in \Ro$ is a common test element for $R$ and $S$, and
$S/\m S$ is $F$-injective, or
\item  $c \in S - \m S$ is a test element for $S$ and $S/\m S$ is $F$-rational, or
\item $R$ is excellent and $S/\m S$ is
$F$-rational.
\endroster
Then $S$ is weakly $F$-regular.
\endproclaim

  We note that the Gorenstein assumption on the fiber is essential, even 
if $R$ is regular.  Even weakening the assumption on the fiber to
$\Q$-Gorenstein is not strong enough to give a good theorem, as Singh
\cite{Si}
gives an example of $R \to S$ flat, where $R$ is a discrete valuation domain,
$S/\m S$ is $\Q$-Gorenstein and strongly $F$-regular, yet $S$ is not weakly
$F$-regular!

  We also prove a corresponding result  for strong
$F$-regularity.  

\proclaim {Theorem 3.6}  Let $(R,\m, K) \to (S, \n, L)$ be a flat
map of $F$-finite reduced rings with Gorenstein closed fiber.  Assume
that $R$ is strongly $F$-regular.  If $S/\m S$ is $F$-rational then
$S$ is strongly $F$-regular.
\endproclaim

In order to prove the first of these theorems we investigate how flat
maps $\phi:(R,\m) \to (S,\n)$ with Gorenstein closed fibers
affect tight closure for $I \inc R$ such that $\l(R/I) < \infty$
and $I$ is irreducible in $R$.  In general these results
do {\it not} depend on the relationship of $R/\m  \to S/\n$ (e.g.,
separability or finiteness).  

While not directly relevant to this paper, we note that other
authors have recently investigated tight closure properties
under good flat maps.  For instance Enescu \cite{En} and Hashimoto
\cite{Ha} have recently shown that for a flat map with $F$-rational
base and $F$-rational closed fiber, the target is $F$-rational
(in the presence of a common test element).

\head \S2.  Background for tight closure \endhead

  Let $R$ be a Noetherian
ring of characteristic $p >0$.  We use $q = p^e$ for a varying power
of $p$ and for an ideal $I\inc R$ we let $I\brq = (i^q:i\in I)$.
Also let $\Ro$ be the complement in $R$ of the union of the
minimal primes of $R$.  Then $x$ is in the \defn{tight closure}
of $I$ if and only if there exists $c \in \Ro$ such that
$cx^q \in I\brq$ for all $q \gg 0$.  If $I^*= I$ then $I$ is
said to be \defn{tightly closed}.  We will say that $I$ is
\defn{Frobenius closed} if $x^q \in I\brq$ for some $q$
always implies that
$x \in I$. 

There is a tight closure
operation for a submodule $N \inc M$, but we will not discuss
this case in general.  It is however useful to discuss tight 
closure in the case of a particular type of direct limit.
Suppose that $M = \varinjlim_t R/I_t$ for a sequence of ideals
$\{I_t\}$.  Let $u \in M$ be an element which is given by
$\{u_t\}$ where in the direct limit system $u_t \mapsto u_{t+1}$.
We will say that $u \in 0^*_M$ if there exists $c \in \Ro$
and a sequence $t_q$ such that for all $q \gg 0$,
$c u_{t_q}^q \in I_{t_q}\brq$.  We will say that $u$ is in the
\defn{finitistic tight closure} of $0$ in $M$, $0^{*fg}_M$, if
there exists $c \in \Ro$ and $t >0$ such that $c u_t^q
\in I_t\brq$ for all $q$.  This definition of finitistic tight
closure agrees with that in
\cite{HH2} for this case.  Clearly $0^{*fg}_M \inc 0^*_M$.

  A ring $R$ in which every ideal is tightly closed is called
\defn{weakly $F$-regular}.  
If every localization of $R$ is
weakly $F$-regular then $R$ is \defn{$F$-regular}.  
When $R$ is reduced then $R^{1/p}$ denotes the
ring of $p$th roots of elements of $R$.  More generally, $R^{1/q}$ is
the ring of $q$th roots.  Clearly $R \inc R^{1/q}$.
If $R$ is
$F$-finite and reduced ($R^{1/p}$ is a finite $R$-module)
then $R$ is called \defn{strongly $F$-regular} if for all
$c \in \Ro$, there exists a $q$ such that
the inclusion $Rc^{1/q} \inc R^{1/q}$ splits over
$R$.  If $R$ is $F$-finite and $R_c$ is strongly $F$-regular for
some $c \in \Ro$, then $R$ is strongly $F$-regular if and only if
there exists $q$ such that $Rc^{1/q} \inc R^{1/q}$ splits over
$R$ \cite{HH1, Theorem 3.3}.  
Strongly $F$-regular rings are $F$-regular, and weakly
$F$-regular rings are normal and under mild conditions (e.g.,
excellent) are Cohen-Macaulay.  

The equivalence of the three
conditions is an important open question.  Let $(R,\m)$ be
an excellent reduced local ring and let $E$ be an injective
hull of the residue field of $R$.  Then $E$ can be written
as a direct limit of the form above since $R$ is approximately
Gorenstein.  Weak $F$-regularity of $R$ is equivalent to
$0_E^{*fg} = 0$ \cite{HH2, Theorem 8.23},
 while strong $F$-regularity is equivalent
to ($F$-finiteness and) $0^*_E = 0$ \cite{LS, Proposition 2.9}. 

By a \defn{parameter ideal} in $(R,\m)$ we mean an ideal generated
by part of a system of parameters.  We say that $(R,\m)$ is
\defn{$F$-rational} if every parameter ideal is tightly closed,
and \defn{$F$-injective} if every parameter ideal is Frobenius
closed  (this is a slightly different notion of $F$-injectivity from
that in \cite{FW}, but is equivalent for CM rings). 
 $F$-rational rings are normal and under mild conditions
are Cohen-Macaulay.   In a Gorenstein ring, $F$-rationality
is equivalent to all forms of $F$-regularity.  

A priori, the multiplier element $c$ in the definition of
tight closure depends on both $I$ and $x$.  If $c$ works for
every tight closure test then we say that $c$ is a \defn{test
element} for $R$.  If $c$ works for every tight closure test
for every completion of every localization of $R$ then we
say that $c$ is a \defn{completely stable test element}.  It
is shown in \cite{HH4} that if $(R,\m)$ is a reduced excellent
domain, $c \in \Ro$, and $R_c$ is Gorenstein and weakly
$F$-regular then $c$ has a power which is a completely stable
test element for $R$.

In \cite{HH2, HH3} it is shown that the multiplier $c$
in the definition of tight closure need not remain constant.
Let $R$ be a domain.
One may have a sequence of elements $c_q$ such
that $c_q x^q \in I\brq$ where $c_q$ must have  ``small order.''
We can obtain a notion of order, denoted ord, by taking
a $\Z$-valued valuation on $R$ which is non-negative on
$R$ and positive on $\m$.  Let $R^+$ be the integral
closure of $R$ in an algebraic closure of the fraction field of
$R$ ($R^+$ has many wonderful properties, such as being a big
Cohen-Macaulay algebra for $R$ when $R$ is excellent \cite{HH5}).
The valuation then extends to 
a function on $R^+$ which takes values in $\Q$.  In particular,
$\ord(c^{1/q}) = \ord(c)/q$.  We will need to use the following
theorem \cite{HH3, Theorem 3.1}.  
\proclaim{Theorem 2.1}  Let $(R,\m)$ be a complete local domain
of characteristic $p$, let $x \in R$ and let $I \inc R$.  Then
$x \in I^*$ if and only if there exists a sequence of elements
$\eps_n \in (R^+)^\circ$ such that $\ord(\eps_n) \to 0$ as $n \to \8$
and $\eps_nx \in IR^+$.
\endproclaim
In fact we would like to strengthen this theorem in order to
apply it to tight closure calculations for non finitely generated
modules which are defined by a direct limit system of ideals.
The proof we give is just an altered version of the proof of
Theorem 3.1 given in \cite{HH3}.  The key component is \cite{HH3,
Theorem 3.3}:
\proclaim{Theorem 2.2}  Let $(R,\m,k)$ be a complete local domain.
Let $\ord$ be a $\Q$-valued valuation on $R^+$ nonnegative on $R$
(and hence on $R^+$) and positive on $\m$ (and, hence, on $\m^+$).
Then there exists a fixed real number $\nu >0$ and a fixed positive
integer $r$ such that for every element $u$ of $R^+$ of order
$< \nu$ there is an $R$-linear map $\phi:R^+ \to R$ such that
$\phi(u) \notin \m^r$.
\endproclaim
The generalization of Theorem 2.1 is given below.
\proclaim{Theorem 2.3}  Let $(R,m)$ be a complete local domain
of characteristic $p$.  Let $M = \varinjlim_t R/I_t$ be an $R$-module
and let $x \in M$.  Suppose that $x$ comes from the sequence 
$\{x_t\}$ where $x_t \mapsto x_{t+1}$.  Then $x \in 0^*_M$ if
and only if there exists a sequence of elements $\eps_n \in
(R^+)^0$ such that $\ord(\eps_n) \to 0$ as $n \to \8$ and
for each $n$ there exists $t$ such that $\eps_n x_t \in I_t R^+$.
\endproclaim
\demo{Proof}
The ``only if'' part is trivial, as if $cx^q = 0$ for all $q \gg 0$ then
we can take $\eps_q = c^{1/q}$.

To see the ``if'' direction,
choose $\nu >0$ and $r$ as in Theorem 2.2.  Fix $q = p^e >0$.
Choose $n$ large enough that $\ord(\eps_n) < \nu/q$.  Let $\eps
= \eps_n^q$.  Then there exists $t$ such that $\eps x_t^q \in I_t\brq
R^+$ and $\ord(\eps) < \nu$.  Applying an $R$ linear map $\phi$
as in Theorem 2.2 we find that $c_q x_t^q \in I_t\brq \inc (I_t\brq)^*$
with $c_q = \phi(\eps) \in R- \m^r$.  Thus, setting
$J_q = \cup_t (I_t\brq)^*:_R x_t^q$ we have $c_q \in J_q$ for
all $q$.

The sequence $J_q$ is nonincreasing.  If for some $t$, $yx_t^{pq}
\in (I_t^{[pq]})^*$ then $c'(y x_t^{pq})^{q'} \in (I_t^{[pq]})^{[q']}
= (I_t^{[pqq']})$ for all $q' \gg 0$ where $c' \ne 0$.  But then
$c'(y x_t^q)^{pq'} \in (I_t \brq)^{[pq']}$ for all $q' \gg 0$ and
hence $y x_t^q \in (I_t \brq)^*$, as required.

Since the sequence $\{J_q\}_q$ is nonincreasing, it cannot have intersection
$0$, or Chevalley's theorem would give $J_q \inc \m^r$ for $q \gg 0$.
As $c_q \in J_q - \m^r$ for all $q$, we can choose a nonzero
element $d \in \cap_q J_q$.  Then for each $q$ there exists $t$ such
that $d x_t^q \in (I_t\brq)^*$.  If $c$ is a test element for $R$
then $cd x_t^q \in I_t\brq$.  Thus $x \in 0^*_M$.
\qed
\enddemo

\proclaim{Proposition 2.4}  Let $(R,\m)$ be an excellent local
domain such that its completion is a domain.  Let $M = \varinjlim_t
R/I_t$ be a direct limit system.  Fix $u \notin 0^*_M$.  Then there
exists $q_0$ such that   $J_q =
\cup_q (I_t\brq:u_t^q)
\inc \m^{[q/q_0]}$ for all $q \gg 0$ (where $\{u_t\}$ represents
$u \in M$ and $u_t \mapsto u_{t+1}$).   In particular if $I \inc R$
we may take $M = R/I$ where the limit system consists of equalities.
Then $u \notin I^*$ implies that $(I\brq:u^q) \inc \m^{[q/q_0]}$.
\endproclaim
\demo{Proof}
Suppose that we can show that the proposition holds in $\hatR$.
Then  $(I_t\brq:_R u_t^q) \inc (I_t\brq:_{\Rhat} u_t^q) \cap R
\inc \m^{[q/q_0]}\Rhat \cap R \inc \m^{[q/q_0]} R$.  Thus
we may assume that $R$ is complete.

For $x \in R$ let $f(x)$ be the largest power of $\m$ that
$x$ is in, and set $\bdf(x) = \lim_{n \to \infty} f(x^n)/n$.  By the valuation
theorem \cite{Re, Theorem 4.16}, there exist a finite number of
$\Z$-valued valuations $v_1,\ldots, v_k$ on $R$ which are non-negative
on $R$ and positive on $\m$ and positive rational numbers
$\vec e k$ such that  $\bdf(x) = \min\{v_i(x)/e_i\}$.  Furthermore,
since $R$ is analytically unramified, there exists a constant $L$ such
that for all $x \in R$, $f(x) \le \lfloor \bdf(x) \rfloor \le
f(x) + L$ (\cite{Re, Theorem 5.32 and 4.16}).   

Now, by Theorem  2.3, for each $v_i$ there
exists a positive real number $\al_i$ such that if $c \in (I_t\brq:u_t^q)$
then $v_i(c) \ge \al_i q$.   Combined with the valuation theorem
we see that $\bdf(c) \ge \min\{q\al_i/e_i \}$.  Let $\al = 
\min\{\al_i/e_i \}$.  Then $f(c) \ge \al q - L -1$.  Let
$s = \mu(\m)$.   Choose $q_1 > 1/\al$, $q_2 \ge L+1$, and $q_3 \ge s$
(all powers of $p$).  Set $q_0 = q_1 q_2 q_3$.  
Then $f(c) \ge \al q - (L+1) \ge q/q_1 - (L+1) \ge q/q_1q_2 - 1
\ge (q/q_0)s -1$.  A simple combinatorial argument shows that
$\m^{(q/q_0)s -1} \in \m^{[q/q_0]}$.  Hence $c \in \m^{[q/q_0]}$.
\qed 
\enddemo 

\head \S3.  Tight closure in flat extension maps \endhead

We show in this section that extending a weakly (respectively, strongly)
$F$-regular ring by a flat map with sufficiently nice Gorenstein
closed fiber yields another weakly (resp., strongly) $F$-regular ring.
These results are Theorems 3.4 and 3.6 (see also Corollary 3.5 for
the $F$-regular case).  

By saying that $\phi: (R,\m) \to (S,\n)$ is flat we mean that $\phi$
is flat and that $\phi(\m) \inc \n$.  Since the map is flat we
then know that given ideals $A,B \inc R$ we have $AS:_S BS = (A:_R B)S$
($B$ finitely generated).  The next lemma merely asserts that
modding out by elements which are regular in the closed fiber
preserves flatness.

\proclaim {Lemma 3.1}  Let $\phi: (R,\m) \to (S,\n)$ be a flat
map.  Let $\vec z d \in S$ be elements whose images in $S/\m S$
are a regular sequence.  Then for any ideal $I$ generated by monomials
in the $z$'s, the ring $S/IS$ is flat over $R$.
\endproclaim

\demo{Proof}  See, for example \cite{HH4, Theorem 7.10a,b}. \qed
\enddemo
The next proposition shows that tight closure behaves well for
irreducible $\m$-primary ideals when extending to $S$.
Given a sequence of elements $\bz = \vec z d$ we will use
$\bz^{[t]}$ to denote $\vec z^t d$.

\proclaim{Proposition 3.2}  Let $\phi:(R, \m,K) \to (S,\n,L)$ be a flat
map with Gorenstein closed fiber.  Let $\bz = \vec z d \in S$ be
elements whose images form a s.o.p.~in $S/\m S$.  Let $I\inc R$
be such that $\l(R/I) < \infty$ and $\dim_K (0:_{R/I} m) = 1$.
Suppose that either 
\roster
\item $R$ and $S$ have a common test element and $S/\m S$ is $F$-injective,
or
\item $c \in S - \m S$ is a test element for $S$, and $S/\m S$ is $F$-rational,
or
\item $R$ is excellent, $\hatR$ is a domain, and $S/\m S$ is $F$-rational.
\endroster
Then  $I$ is tightly closed in $R$ $\iff$ for all $t >0$,
$IS + (\bz)\brt S$ is tightly closed in $S$ $\iff$ there exists $t >0$
such that $IS + (\bz)\brt S$ is tightly closed in $S$.
\endproclaim

\demo{Proof}
Let $b \in S$ have as its image the socle element in $S/\m S + (\bz)S$.
Let $u \in R$ be the socle element mod $I$.  Then the socle element
of $S/(IS + (\bz)S)$ is $ub$ since the map $R/I \to R/I \otimes S = S/IS$
is flat with Gorenstein fibers (there is only one fiber).  

Suppose that $I$ is tightly closed.  There is no loss of generality
in taking $t=1$.  If $IS + (\bz)S$ is not tightly closed in $S$ then
we have $c(ub)^q \in (I\brq + (\bz)\brq)S$ for all $q$.  
In case (1) we may take $c \in \Ro$, so that  
$$
b^q \in (I\brq + (\bz)\brq)S:_S cu^q = (I\brq:_R cu^q)S + (\bz)\brq S
\inc \m S + (\bz)\brq S
$$
for all $q \gg 0$.  The first equality is a consequence of flatness,
while the inclusion follows since $u \notin I^*$.  By our assumption
that $S/\m S$ is $F$-injective we reach the contradictory conclusion
that $b \in ((\bz) + \m)S$.
In case (2) we have
$$
cb^q \in (I\brq + (\bz)\brq)S:_S u^q = (I\brq:_R u^q)S + (\bz)\brq S
\inc \m S + (\bz)\brq S
$$
for all $q \gg 0$.  As $S/\m S$ is $F$-rational, it is a domain, so
$c\ne 0$ in $S/\m S$.  This contradicts our hypothesis that $S/\m S$ is
$F$-rational (in fact it is enough to assume that $I$ is Frobenius
closed to reach this conclusion).  In case (3) we can choose $q_0$
as in Proposition 2.4, and then
$$
c(b^{q_0})^{q/q_0} \in (I\brq + (\bz)\brq)S:_S u^q = 
(I\brq:_R u^q)S + (\bz)\brq S \inc \m^{[q/q_0]}S + ((\bz)^{[q_0]})^{[q/q_0]}
$$ for all $q/q_0$.  But then $b^{q_0} \in (\m S + (\bz)^{[q_0]})^*$.
By persistence, the image of $b^{q_0}$ is in  $((\bz)^{[q_0]}S/\m S)^*$,
which contradicts the $F$-rationality of $S/\m S$.

Suppose now that $IS + (\bz)\brt S$ is tightly closed in $S$ for
all $t$, but $I$ is
not tightly closed in $R$.  Then $u \in (IR)^* \inc (I + (\bz)\brt)^*$
(since $\Ro \inc S^0$).  But then $u \in \cap_t (IS + (\bz)\brt S)^* \cap R
\inc \cap_t (IS + (\bz)\brt S) \cap R \inc IS \cap R = IR$.

Finally, suppose that $(IS + (\bz)^{[t_0]})S$ is tightly closed
for some $t_0$.  Given any $t$, the socle element of 
$(IS + (\bz)\brt)S$ is $(z_1\cdots z_d)^{t-1}ub$.  If
$c ((z_1\cdots z_d)^{t-1}ub)^q \in (IS + (\bz)\brt)\brq$ then
by flatness, $c ((z_1\cdots z_d)^{t_0-1}ub)^q \in (IS + (\bz)^{[t_0]})\brq$.
Therefore, one such ideal tightly closed shows that all such
ideals are tightly closed.
\qed
\enddemo

To deal with strong $F$-regularity we need to give a similar proposition
with $R/I$ replaced by the injective hull $E_R(R/\m)$.  Suppose
that we can write $E = E_R(R/\m) = \varinjlim_t R/J_t$, the set $\{u_t\} \inc R$
is a collection of elements such that $u_t \mapsto u_{t+1}$ in
the map $R/J_t \to R/J_{t+1}$ and the image of each $u_t$ in $E$
is the socle element of $E$.  It suffices that $R$ be approximately
Gorenstein \cite{Ho2} (e.g., excellent and normal, or even reduced)
 to obtain $E$ in this manner.  In particular an $F$-finite ring
is excellent \cite{Ku}, so a reduced $F$-finite ring is approximately
Gorenstein.

\proclaim{Proposition 3.3}  Let $(R,\m, K) \to (S,\n, L)$ be a flat
map of $F$-finite reduced rings with Gorenstein closed fiber.  
\roster
\item   If $R c^{1/q} \inc \Rq$ splits for some $q$ (over $R$) and $S/\m S$
is $F$-injective then $S c^{1/q} \inc \Sq$ splits for some $q$ (over $S$).
\item  If $0$ is Frobenius closed in $E_R(K)$, $S/\m S$ is $F$-rational
and $c \in S - \m S$ then there exists $q$ such that $S c^{1/q} \inc
\Sq$ splits (over $S$).
\endroster
\endproclaim

\demo{Proof}
Choose $\bz = \vec z d \in S$ elements which generate a s.o.p.~in 
$S/\m S$.  By \cite{HH4, Lemma 7.10} we have $E_S(L) =
\varinjlim_{v} S/(\bz^{[v]}) \otimes_R  E_R(K) = \varinjlim_{v,t}
S/(\bz^{[v]}) \otimes_R R/J_t  = \varinjlim_t  S/(\bz\brt,J_t)S$.
If $b \in S$ generates the socle element in $S/(\m +(\bz))S$ then the
image of $(z_1\cdots z_d)^{t-1}bu_t$ in $S/((\bz\brt) + J_t)S$
maps to the socle element of $E_S$ (where $u_t$ is as given above).

In case (1), if for all $q$ the inclusion $Sc^{1/q} \to \Sq$ fails
to split, by \cite{Ho1, Theorem 1 and Remark 2} 
for all $q$ there exists $t_q$ such that
$$
c(z_1 \cdots z_d)^{(t_q-1)q} b^q u_{t_q}^q \in ((\bz^{[t_q]}),J_{t_q})\brq
S.
$$
Hence $(z_1 \cdots z_d)^{(t_q-1)q} b^q \in 
((\bz),J_{t_q})\brq :_S c u_{t_q}^q
\inc (J_{t_q}\brq:_R c u_{t_q}^q)S + (\bz^{[t_q]})\brq S
\inc \m S + (\bz^{[t_q]})\brq S$ for $q \gg 0$ (we are using
here that if $Rc^{1/q} \inc \Rq$ splits for some $q$ then 
$R c^{1/q'} \inc R^{1/q'}$ splits for all $q' \ge q$).
Thus $b^q \in \m S + (\bz)\brq$ since $S/\m S$ is CM.  This contradicts
the $F$-injectivity of $S/ \m S$.

To see (2), if there is no splitting we obtain
$$ c(z_1 \cdots z_d)^{(t_q-1)q} b^q \in 
(\bz^{[t_q]},J_{t_q})\brq :_S  u_{t_q}^q
\inc (J_{t_q}\brq:_R  u_{t_q}^q)S + ((\bz^{[t_q]})\brq S
\inc \m S + ((\bz^{[t_q]})\brq S
$$
 and hence $cb^q \in \m S + (\bz)\brq$.  This contradicts the
$F$-rationality of $S/\m S$.
\qed
\enddemo

We can now give our main theorems on the extension of weakly and strongly
$F$-regular rings by flat maps with Gorenstein closed fiber.

\proclaim{Theorem 3.4}  Let $\phi:(R,\m) \to (S,\n)$ be a flat map.  
Assume that
$S/\m S$ is Gorenstein and $R$ is weakly $F$-regular and CM.  
Suppose that either
\roster
\item  $c \in \Ro$ is a common test element for $R$ and $S$, and
$S/\m S$ is $F$-injective, or
\item  $c \in S - \m S$ is a test element for $S$ and $S/\m S$ is $F$-rational, or
\item $R$ is excellent and $S/\m S$ is
$F$-rational.
\endroster
Then $S$ is weakly $F$-regular.
\endproclaim

\demo{Proof}  To see that $S$ is weakly $F$-regular it suffices to show
that there exists a sequence of irreducible tightly closed ideals
of $S$ cofinite with the powers of $\n$.  As $R$ is weakly $F$-regular
(so normal) and CM it is approximately Gorenstein.  Say that $\{J_t\}$ is a
sequence of irreducible ideals cofinite with the powers of $\m$.
Let $\bz = \vec z d \in S$  be elements which form a s.o.p.~in $S/\m S$.
Then $(J_t + \bz\brt)S$ is a sequence of irreducible ideals in $S$
cofinal with the powers of $\n$.  By Proposition 3.2, in
cases (1), (2), and (3), the ideals $(J_t + \bz\brt)S$ are tightly
closed in $S$ (in case (3), $\hatR$ is still weakly $F$-regular, so
is a domain).  Therefore $S$ is weakly $F$-regular. 
We note that in case (2) we may weaken the assumption that $R$ is
weakly $F$-regular to the assumption that $R$ is $F$-pure (see the comment 
in the proof of Proposition 3.2, part (2)). \qed
\enddemo

The next corollary should be compared with \cite{HH4, Theorem 7.25(c)}.

\proclaim{Corollary 3.5}  Let $(R,\m) \to (S, \n)$ be a flat map
of excellent rings
with Gorenstein  fibers.  Suppose that the generic fiber is
$F$-rational and all other fibers are $F$-injective.
  If $R$ is $F$-regular  then
$S$ is $F$-regular.
\endproclaim

\demo{Proof}
By hypothesis the generic fiber is Gorenstein and $F$-rational, therefore
there is a $c \in \Ro$ which is a common completely stable 
test element.  $F$-regularity
is local on the prime ideals of $S$ and the fiber of such a localization
is the localization of a fiber, hence Gorenstein and $F$-injective (the
property of $F$-injectivity is easily seen to localize).
Therefore Theorem 3.4(1) always applies.  \qed
\enddemo

\proclaim {Theorem 3.6}  Let $(R,\m, K) \to (S, \n, L)$ be a flat
map of $F$-finite reduced rings with Gorenstein closed fiber.  Assume
that $R$ is strongly $F$-regular.  If $S/\m S$ is $F$-rational then
$S$ is strongly $F$-regular.
\endproclaim

\demo{Proof} We must show that there exists an element $c \in S^0$ such that
$S_c$ is strongly $F$-regular and $Sc^{1/q} \inc \Sq$ splits for some
$q$.  

If there exists $c \in \Ro$ such that $S_c$ is strongly $F$-regular
(i.e., a power of $c$ is a common test element for $R$ and $S$)
then we are done by Proposition 3.3(1).  Even if $R$ and $S$ have no (apparent)
common
test element, however,  we claim that there exists $c \in S - \m S$ such that
$S_c$ is strongly $F$-regular.  Once we have shown this, the theorem
follows by Proposition 3.3(2).

Since the non-strongly $F$-regular locus is closed \cite{HH1, Theorem
3.3} it suffices
to show that $S_{\m S}$ is strongly $F$-regular, for then there 
exists an element $ c \in S - \m S$ such that $S_c$ is 
strongly $F$-regular.  Let $B = S_{\m S}$.
Then $R \to B$ is flat and the closed fiber is a field.  In
particular $E_B(B/\m B) = E_R(K) \otimes_R B$.  As $R$ is
strongly $F$-regular (so normal) it is approximately Gorenstein.
Say $E_R
= \varinjlim_t R/J_t$ with socle element mapped to by $u_t$ (as before).
Then $u_t \in B/J_t B$ will still map to the socle element $u$ in 
$E_B$.   Suppose that $u \in 0^*_{E_B}$.
  This means there exists $b \in B_0$ such that for all $q$ there exists $t_q$
such that  $b u_{t_q}^q \in J_{t_q}\brq B$.  Hence
$b \in J_{t_q}\brq:_B u_{t_q}^q = (J_{t_q}\brq:_R u_{t_q}^q)B$.
Note that $R$ is an excellent normal domain, so its completion
remains a domain.  Thus by Proposition 2.4 we see that
 as $q \to \infty$, $(J_{t_q}\brq:_R u_{t_q}^q)$ gets
into larger and larger powers of the maximal ideal, since
$0$ is tightly closed in $E_R$.  Thus $b \in \cap_N \m^N B = 0$,
a contradiction.  \qed
\enddemo

\enddocument